\documentstyle[12pt,epsf]{article}
\newcommand{\rr}{{\rm{I \! R}}}
\newcommand{\rn}{{\rm{I \! N}}}

\newcommand{\ds}{\displaystyle}
\newcommand{\f}{\frac}
\newcommand{\lam}{\lambda}

\newcommand{\om}{\omega}

\newcommand{\ca}{{\cal A}}
\newcommand{\th}{\theta}
\newcommand{\cc}{{\rm{{\footnotesize{l}}\!\!\! C}}}

\title{\bf A dynamic p53-mdm2 model with delay kernel }
\author{\small{R.F. HORHAT$^a$\thanks{Corresponding author}, M. NEAM\c TU$^{b}$, D. OPRI\c S$^{c}$}}
\date{}

\begin{document}
\maketitle

\begin{tabular}{cccccccc}
\scriptsize{$^{a}$ Department of Biophysics and Medical Informatics,}\\
\scriptsize{University of Medicine and Pharmacy, Piata Eftimie Murgu, nr. 3, 300041, Timi\c soara, Romania,}\\
\scriptsize{E-mail: rhorhat@medinfo.umft.ro}\\
\scriptsize{$^{b}$Department of Economic Informatics and Statistics, Faculty of Economics,}\\
\scriptsize{West University of Timi\c soara, str. Pestalozzi, nr. 16A, 300115, Timi\c soara, Romania,}\\
\scriptsize{E-mail:gabriela.mircea@fse.uvt.ro, mihaela.neamtu@fse.uvt.ro,}\\
\scriptsize{$^{c}$ Department of Applied
Mathematics, Faculty of Mathematics,}\\
\scriptsize{West University of Timi\c soara, Bd. V. Parvan, nr. 4, 300223, Timi\c soara, Romania,}\\
\scriptsize{E-mail: opris@math.uvt.ro}\\

\end{tabular}

\bigskip

\begin{center}
\small{{\bf Abstract}. }
\end{center}

\begin{quote} \small{ Specific activator and repressor
transcription factors which bind to specific regulator DNA
sequences, play an important role in gene activity control.
Interactions between genes coding such transcripion factors should
explain the different stable or sometimes oscillatory gene
activities characteristic for different tissues. In this paper,
the dynamic P53-Mdm2 interaction model with distributed delays and
weak kernel, is investigated. Choosing the delay or the kernel's
coefficient as a bifurcation parameter, we study the direction and
stability of the bifurcating periodic solutions. Some numerical
examples are finally given for justifying the theoretical
results.}
\end{quote}

\noindent{\small{{\it Keywords:} delay differential equation,
stability, Hopf bifurcation, P53, MDM2.

}}

\section*{\normalsize\bf 1. Introduction}

\hspace{0.6cm}

P53 is a very important gene in oncogenesis. It is also known as
"Guardian of the genome". Its anomalies are almost universal in
tumoral cells [5]. The full activity of p53 gene starts when is
detected a DNA damage [3,4]. These damages are mainly formed by
DSB (Double-Strand Break) lesions [9]. Around these DSBs it will
be formed repair complexes. These complexes include at eukaryotes
the proteins: Mre11, Rad50 and NBS1 (MRN complex) and they are the
signal for activation of ATM. The DSBs repair protein complexes
count only for the initial activation of ATM, because the main
activation is an autocatalytical process. ATM, at its own,
represents the signal for activation of gene p53. Depending on the
level of ATM, p53 will lead to two outcomes for the cell: one is
the cell cycle arrest induced by a low level or a brief elevation
of p53 protein, and the other is the apoptosis induced by a high
level or a prolonged elevation of p53 protein [8]. Each of these
two outcomes could not be an option for all the cells. For
example, apoptosis is not accepted for neurons or myocardial
muskular cells because they do not divide in adult life, so these
cells will choose for cell cycle arrest. On the other hand for the
enterocytes (cells of digestive tube) apoptosis is a common
option, because these cells divide themselves very quickly and
their lifetime is no longer than 2 days. Now is clear that should
be a very good control of p53 activity in such a manner that the
cell goes on right pathway (i.e. apoptosis or cell cycle arrest).
This control is achieved with the help of mdm2 gene with which p53
makes a feedback loop [8, 12].

In the last years, the approaches of P53 dynamics as response to
DNA damage comprise modelings in which are described three
distinct subsystems: a DNA damage repair module, an ataxia
telengiectasia mutated (ATM) switch and the P53-Mdm2 oscillator.

In what follows we will consider a model only for the third
module. The variables of the model are: $x_1$ P53-mRNA
concentration, $x_2$ Mdm2-mRNA concentration, $y_1$ P53-protein
concentration and $y_2$ Mdm2-protein concentration.

We consider P53-Mdm2 model with kernel delay given by:

$$\begin{array}{l}
\vspace{0.1cm}
\dot x{}_1(t)=a_1-a_2x_1(t),\\
\vspace{0.1cm}
\dot y{}_1(t)=b_1x_1(t)-b_2y_1(t)-b_{12}y_1(t)y_2(t),\\
\dot x{}_2(t)=\int^{\infty}_0 k_1(s)f(y_1(t-s))ds-c_2x_2(t),\\
\dot
y{}_2(t)=\int_0^{\infty}k_2(s)x_2(t-s)ds-d_2y_2(t)-d_{12}y_1(t)y_2(t)\end{array}\eqno(1)$$
where:  $a_2, c_2$ are the rates for mRNA degradation, $b_2, d_2,
b_{12}, d_{12}$ are the rates for proteins degradation. The
function $f:\rr_+\rightarrow\rr_+$, is the Hill function, given
by:
$$f(x)=\f{x^n}{a^n+x^n}$$
with $n\in\rn^*, a>0.$ The parameters $a_1$, $a_2$, $b_1$, $b_2$,
$c_2$, $b_{12}$, $d_2$ $d_{12}$ of the model are assumed to be
positive numbers less or equal to one, the memory functions $k_1$,
$k_2$ that reflect the influence of the past states on the current
dynamics are a nonnegative bounded function defined on
$[0,\infty)$ and
$$\int_0^{\infty}k_i(s)ds=1, \quad \int_0^{\infty}sk_i(s)ds<\infty, \quad i=1,2.$$ The memory function is called delay
kernel. The delay becomes a discrete one when the delay kernel is
a delta function at a certain time. Usually, we employ the
following form:
$$k_i(s)=\ds\f{q^{p+1}_i}{p!}s^pe^{-q_is}, \quad q_i>0, i=1,2, p\geq0$$ for the memory
function. When $p=0$ and $p=1$ the memory function are called
"weak" and "strong" kernel respectively.

For $k_1(s)=\delta(s-\tau_1)$, $k_2(s)=\delta(s-\tau_2)$,
$\tau_1\geq0$, $\tau_2\geq0$ the model is given by:

$$\begin{array}{l}
\vspace{0.1cm}
\dot x{}_1(t)=a_1-a_2x_1(t),\\
\vspace{0.1cm}
\dot y{}_1(t)=b_1x_1(t)-b_2y_1(t)-b_{12}y_1(t)y_2(t),\\
\dot x{}_2(t)=f(y_1(t-\tau_1))-c_2x_2(t),\\
\dot
y{}_2(t)=x_2(t-\tau_2)-d_2y_2(t)-d_{12}y_1(t)y_2(t).\end{array}\eqno(2)$$

In (2) with $\tau_1=\tau$, $\tau_2=0$, $d_{12}=0$ we obtain the
model from [11] and it was studied in [10] which suggests that
there is an oscillatory behavior based on observations obtained
using only numerical simulations.

In present paper we will analyze the model (1) with $d_{12}=0$
with the following initial values:
$$\begin{array}{l}
x_1(0)=\bar x_1, y_1(\th)=\varphi_1(\th), \th\in(-\infty,0],\\
x_2(\theta)=\varphi_2(\theta), \th\in(-\infty,0], y_2(0)=\bar
y_2,\end{array}$$ with $\bar x_1\geq0$,  $\bar y_2\geq0$,
$\varphi_1(\theta)\geq0$, $\varphi_2(\theta)\geq0$, for all
$\theta\in(-\infty, 0]$ and $\varphi_1$, $\varphi_2$ are
differentiable functions.

Also, two delays appear in the leukopoiesis model which is
analyzed in [1].

The paper is organized as follows. In section 2, we discuss the
local stability for the equilibrium state of system (1), with
$d_{12}=0$, if the delay kernels $k_1$, $k_2$ are delta functions,
and $k_1$ is delta function, $k_2$ is weak kernel function. We
investigate the existence of the Hopf bifurcation with respect to
the parameters of the delay kernels $k_1$, $k_2$. In section 3,
the direction of Hopf bifurcation is analyzed by the normal form
theory and the center manifold theorem introduced by Hassard [4].
Numerical simulations for justifying the theoretical results are
illustrated in section 4. Finally, some conclusions are made.

\section*{\normalsize\bf{2. Local stability and the existence of the Hopf bifurcation.}}

\hspace{0.6cm} We consider the model:

$$\begin{array}{l}
\vspace{0.1cm}
\dot x{}_1(t)=a_1-a_2x_1(t),\\
\vspace{0.1cm}
\dot y{}_1(t)=b_1x_1(t)-b_2y_1(t)-b_{12}y_1(t)y_2(t),\\
\dot x{}_2(t)=\int^{\infty}_0 k_1(s)f(y_1(t-s))ds-c_2x_2(t),\\
\dot
y{}_2(t)=\int_0^{\infty}k_2(s)x_2(t-s)ds-d_2y_2(t).\end{array}\eqno(3)$$

{\bf Proposition 1}. {\it If $b_2^2<b_1$ and
$y_{10}\in(0,\ds\f{a_1b_1}{a_2b_2})$ is a solution of equation

$$\alpha x^{n+1}-\beta x^n+\gamma x-\delta=0$$
where $$\alpha=a_2(b_{12}+b_2c_2d_2), \beta=a_1b_1c_2d_2,
\gamma=a_2b_2c_2d_2a^n, \delta=a_1b_1c_2d_2a^n$$ then the
equilibrium point $X^\ast$ of system (3) has the coordinates:
$$x_{10}=\ds\f{a_1}{a_2}, \quad x_{20}=d_2y_{20},
\quad y_{20}=\ds\f{a_1b_1-a_2b_2y_{10}}{b_{12}y_{10}a_2}.$$}

We consider the following translation:
$$x_1=u_1+x_{10}, y_1=u_2+y_{10}, x_2=u_3+x_{20},
y_2=u_4+y_{20}.\eqno(4)$$ With respect to (4), system (3) can be
expressed as:
$$\begin{array}{l}
\vspace{0.1cm}
\dot u{}_1(t)=-a_2u_1(t),\\
\vspace{0.1cm}
\dot u{}_2(t)=b_1u_1(t)-(b_2+b_{12}y_{20})u_2(t)-b_{12}y_{10}u_4(t)-b_{12}u_2(t)u_4(t),\\
\dot u{}_3(t)=\int_0^\infty k_1(s) f(u_2(t-s)+y_{10})ds -c_2(u_3(t)+x_{20}),\\
\dot u{}_4(t)=\int_0^\infty
k_2(s)u_3(t-s)ds-d_2u_4(t).\end{array}\eqno(5)$$

System (5) has $(0,0,0,0)$ as equilibrium point.

To investigate the local stability of the equilibrium state we
linearize system (5). We expand it in a Taylor series around the
origin and neglect the terms of higher order than the first order
for the functions from the right side of (5). We obtain:
$$\dot U(t)=AU(t)+B_1U_1(t)+B_2U_2(t),\eqno(6)$$
where
$$A\!\!=\!\!\left(\!\!\!\!\begin{array}{cccc}
\vspace{0.2cm}
-a_2 & 0 & 0 & 0\\
\vspace{0.2cm}
b_1 & -(b_2\!\!+\!\!b_{12}y_{20}) & 0 & -b_{12}y_{10}\\
0 & 0 & -c_2 & 0\\
\vspace{0.2cm} 0 & 0 & 0 & -d_2
\end{array}\right)\eqno(7)$$
$$B_1=\left(\begin{array}{cccc} \vspace{0.2cm}
0 & 0 & 0 & 0\\
\vspace{0.2cm}
0 & 0 & 0 & 0\\
\vspace{0.2cm}
0 & \rho & 0 & 0\\
\vspace{0.2cm} 0 & 0 & 0 & 0\end{array}\right), \quad
B_2=\left(\begin{array}{cccc} \vspace{0.2cm}
0 & 0 & 0 & 0\\
\vspace{0.2cm}
0 & 0 & 0 & 0\\
\vspace{0.2cm}
0 & 0 & 0 & 0\\
\vspace{0.2cm} 0 & 0 & 1 & 0\end{array}\right) \eqno(8)$$ with
$\rho=f'(y_{10})$, $U(t)=(u_1(t)$,$u_2(t)$,$u_3(t)$,$u_4(t))^T$,
\vspace{0.2cm}
$U_i(t)=(\int_0^{\infty}\rho_i(s)u_1(t-s)ds$,$\int_0^{\infty}\rho_i(s)u_2(t-s)ds$,$\int_0^{\infty}\rho_i(s)u_3(t-s)ds$,
$\int_0^{\infty}\rho_i(s)u_4(t-s)ds)^T, i=1,2.$
\smallskip

The characteristic equation corresponding to system (6) is
$\Delta(\lambda)=0$, where $$\Delta(\lambda)=det(\lambda
I-A-(\int_0^{\infty}k_1(s)e^{-\lambda
s}ds)B_1-(\int_0^{\infty}k_2(s)e^{-\lambda s}ds)B_2).\eqno(9)$$

From (7), (8) and (9) it results:
$$\Delta(\lambda)=(\lambda+a_2)\Delta_1(\lambda)$$ where
$$\Delta_1(\lambda)=\lam^3+p_2\lam^2+p_1\lam+p_0+r(\int _0^\infty k_1(s) e^{-\lam s}ds)(\int _0^\infty k_2(s) e^{-\lam s}ds)\eqno(10)$$
with
$$\begin{array}{l}
p_2=b_2+c_2+d_2+b_{12}y_{20}, p_1=(c_2+d_2)(b_2+b_{12}y_{20})+c_2d_2\\
p_0=c_2d_2(b_2+b_{12}y_{20}), r=\rho
b_{12}y_{10}.\end{array}\eqno(11)$$

The equilibrium point $X^*=(x_{10}$, $y_{10}$, $x_{20}$,
$y_{20})^T$ is locally asymptotically stable if and only if all
eigenvalues of $\Delta(\lambda)=0$ have negative real parts.
Because $a_2>0$, we will analyze the equation
$\Delta_1(\lambda)=0$. The analysis of the sign of real parts of
eigenvalues is complicated and a direct approach cannot be
considered.

We will analyze the eigenvalues for the equation
$\Delta_1(\lambda)=0$ if the delay kernels $k_1$ and $k_2$ are
delta functions or $k_1$ is delta function and $k_2$ is weak
function.

{\bf Proposition 2.} {\it If $k_1(s)=\delta(s-\tau_1)$,
$k_2(s)=\delta(s-\tau_2)$, $\tau_1\geq0$, $\tau_1\geq0$ then:

(i) function (10) is given by:
$$\Delta_1(\lambda, \tau)=\lam^3+p_2\lam^2+p_1\lam+p_0+re^{-\lam
\tau}\eqno(12)$$where $\tau=\tau_1+\tau_2.$

(ii) if $\tau_1=0$, $\tau_2=0$ then the equilibrium state $X^\ast$
of system (5) is locally asymptotically stable if and only if
$$p_1p_2>p_0+r\eqno(13)$$ where $p_1$, $p_2$, $p_0$, r are given
by (11).}

We are looking for the values $\tau_0$ so that the equilibrium
point $X^*$ changes from local asymptotic stability to instability
or vice versa. If the relation (13) holds, for $\tau>0$
sufficiently small all roots of $\Delta_1(\lambda,\tau)=0$ have
negative real parts.The critical delay $\tau_0$ is the smallest
positive value of $\tau$ where $\Delta_1(\lambda,\tau)=0$ has
imaginary roots. Let $\lambda=\pm i\om$ be these solutions with
$\omega>0$. Separating real and imaginary parts of
$\Delta_1(i\om,\tau)=0$ we obtain:

$$rcos(\om\tau)=p_2\om^2-p_0, \quad rsin(\om\tau)=\om p_1-\om^3.\eqno(14)$$

A solution of (14) is a pair $(\om_0,\tau_0)$where $\om_0$ is a
positive root of the equation:

$$x^6+(p_2^2-2p_1)x^4+(p_1^2-2p_0p_2)x^2+p_0^2-r^2=0$$
and $\tau_0$ is given by:
$$\tau_0=\ds\f{1}{\om_0}arctg\ds\f{\om_0(-p_1+\om_0^2)}{-p_2\om_0^2+p_0}.$$

From (12), we obtain:
$$\lambda^{'}=\ds\f{d\lambda}{d\tau}=-\ds\f{\lambda
r}{e^{\lambda\tau}(3\lambda^2+2p_2\lambda+p_1)-r\tau}.\eqno(15)$$

Then, we evaluate (15) at $\lam=i\om_0$ and $\tau=\tau_0$ and
obtain:
$$\lambda^{'}(\tau_0)=\ds\f{\om_0rl_2}{l_1^2+l_2^2}+i\ds\f{\om_0rl_1}{l_1^2+l_2^2},$$
where $$\begin{array}{l}
l_1=(p_1-3\om_0^2)cos(\om_0\tau_0)-2p_2\om_0sin(\om_0\tau_0)-r\tau_0\\
l_2=(p_1-3\om_0^2)sin(\om_0\tau_0)+2p_2\om_0cos(\om_0\tau_0).\end{array}$$

From the above analysis and the standard Hopf bifurcation theory
[4], we have the following result:

{\bf Proposition 3}. {\it If $p_1$, $p_2$, $p_0$, r satisfy (13)
and $p_0<r$, for $\tau=\tau_0$, $\om=\om_0$ then: $$Re \left
(\ds\f{d\lambda}{d\tau}\right)_{\lam=i\om_0,
\tau=\tau_0}=\ds\f{\om_0rl_2}{l_1^2+l_2^2}\neq0.$$ and a Hopf
bifurcation occurs at the equilibrium state $X^{\ast}$ as $\tau$
passes through $\tau_0$.}

{\bf Proposition 4}. {\it If $k_1(s)=\delta(s-\tau_1)$,
$k_2(s)=q_2e^{-sq_2}$, $\tau_1\geq0$, $q_2>0$ then:

(i) function (10) is given by:
$$\Delta_1(\lambda, \tau_1)=(\lam+q_2)(\lam^3+p_2\lam^2+p_1\lam+p_0)+rq_2e^{-\lam
\tau_1};\eqno(16)$$

(ii) if $\tau_1=0$, then the equilibrium state $X^\ast$ of system
(5) is locally asymptotically stable if and only if
$$\begin{array}{l}
D_2=(p_2+q_2)(p_1+q_2p_2)-(p_0+p_1q_2)>0,\\
D_3=(p_0+p_1q_2)D_2-(p_2+q_2)^2(q_2p_0+rq_2)>0.
\end{array}\eqno(17)$$}

We are looking for the values $\tau_{10}^{\ast}$ so that the
equilibrium point $X^*$ changes from local asymptotic stability to
instability or vice versa. The critical delay $\tau_{10}^{\ast}$
is the smallest positive value of $\tau_1$ where
$\Delta_1(\lambda,\tau_{10}^{\ast})=0$ has imaginary roots. Let
$\lambda=\pm i\om$ be these solutions with $\omega>0$. Separating
real and imaginary parts of $\Delta_1(i\om,\tau)=0$ we obtain:

$$\begin{array}{l}
q_2rcos(\om \tau_1)=-\om^4+(p_1+q_2p_2)\om^2-q_2p_0,\\
q_2rsin(\om\tau_1)=-(p_2+q_2)\om^3+(p_0+p_1q_2)\om.\end{array}\eqno(18)$$

A solution of (18) is a pair $(\om_{10},\tau_{10}^{\ast})$where
$\om_{10}$ is a positive root of the equation:

$$x^8+n_1x^6+n_2x^4+n_3x^2+n_4=0$$ where
$$\begin{array}{l}
n_1=(p_2+q_2)^2-2(p_1+q_2p_2),\\
n_2=(p_1+q_2p_2)^2+2q_2p_0-2(p_0+p_1q_2)(p_2+q_2),\\
n_3=(p_0+p_1q_2)^2-2q_2p_0(p_1+q_2p_2),\\
n_4=p_0^2q_2^2-r^2q_2^2.
\end{array}$$ and $\tau_{10}^{\ast}$ is given by
$$\tau_{10}^{\ast}=\ds\f{1}{\om_{10}}arctg\ds\f{(p_2+q_2)\om_{10}^3-(p_0+p_1q_2)\om_{10}}{\om_{10}^4-(p_1+q_2p_2)\om_{10}^2+q_2p_0}.\eqno(19)$$

From (16), we obtain:
$$\lambda^{'}\!\!=\!\!\ds\f{d\lambda}{d\tau_1}\!\!=\!\!\ds\f{\lambda
rq_2}{e^{\lambda\tau_1}(4\lambda^3+3(p_2+q_2)\lam^2+2(p_1+q_2p_2)\lambda+p_0+p_1q_2)-rq_2\tau_1}.\eqno(20)$$

Then, we evaluate at $\lam=i\om_{10}$ and
$\tau_1=\tau_{10}^{\ast}$ and obtain:
$$\lambda^{'}(\tau_{10}^{\ast})=-\ds\f{\om_{10}rq_2l_{20}}{l_{10}^2+l_{20}^2}i+\ds\f{\om_{10}rq_2l_{10}}{l_{10}^2+l_{20}^2},$$
where $$\begin{array}{ll} l_{10}&
=(-3(p_2+q_2)\om_{10}^2+p_0+p_1q_2)cos(\om_{10}\tau_{10})\\
& +(4\om_{10}^3-2(p_1+q_2p_2)
\om_{10})sin(\om_{10}\tau_{10})-r\tau_1q_2\\
l_{20}
&=(-3(p_2+q_2)\om_{10}^2+p_0+p_1q_2)sin(\om_{10}\tau_{10})\\
&+
(4\om_{10}^3+2(p_1+q_2p_2)\om_{10})cos(\om_0\tau_0).\end{array}$$

We have the following result:

{\bf Proposition 5}. {\it If $p_1$, $p_2$, $p_0$, r, $q_2$ satisfy
(17) for $\tau_1=\tau_{10}^{\ast}$, $\om=\om_{10}$ then:
$$Re \left (\ds\f{d\lambda}{d\tau_1}\right)_{\lam=i\om_{10},
\tau_1=\tau_{10}^{\ast}}=\ds\f{\om_{10}rq_2l_{10}}{l_{10}^2+l_{20}^2}\neq0.$$
and a Hopf bifurcation occurs at the equilibrium state $X^{\ast}$
as $\tau_1$ passes through $\tau_{10}^{\ast}$.}

\section*{{\normalsize\bf 3. Direction and stability of the Hopf bifurcation}}

In what follows, we will study the direction and stability in two
cases: in the first case the both kernels are delta function and
in the second case the kernel $k_1$ is delta function and the
kernel $k_2$ is weak function.

\vspace{0.2cm} {\bf 3.1. The case $k_1(s)=\delta(s-\tau_1)$,
$k_2(s)=\delta(s-\tau_2)$, $\tau_1\geq0$, $\tau_2\geq0$.}

\hspace{0.6cm} For $k_1(s)=\delta(s-\tau_1)$,
$k_2(s)=\delta(s-\tau_2)$, $\tau_1\geq0$, $\tau_1\geq0$ from
Proposition 2, we obtained some conditions which guarantee that
system (5) undergoes Hopf bifurcation at $\tau=\tau_0$. In this
section, we study the direction, the stability and the period of
the bifurcating periodic solutions. The used method is based on
the normal form theory and the center manifold theorem introduced
by Hassard [4]. We know that if $\tau=\tau_0$ then all roots of
$\Delta_1(\lambda,\tau_0)=0$, when $\Delta_1(\lambda,\tau_0)$ is
given by (12), other than $\pm i\om_0$ have negative real parts
and any roots of the form $\lambda(\tau)=\alpha(\tau)\pm
i\om(\tau)$ satisfies $\alpha(\tau_0)=0$, $\om(\tau_0)=\om_0$ and
$\ds\f{d\alpha(\tau_0)}{d\tau}\neq0.$

Suppose that for given $a_1$, $a_2$, $b_1$, $b_2$, $b_{12}$,
$c_2$, a, $d_2$ there is $\tau_0$ for which
$\Delta_1(\lambda,\tau_0)=0$ exhibits a Hopf bifurcation. We
consider $\tau_{10}=\tau_0-\tau_2$, where $\tau_2<2\tau_0$ and
$\tau_1=\tau_{10}+\mu$, $\mu\in\rr$. We regard $\mu$ as the
bifurcation parameter.

For $\Phi\in C^1=C([-\tau_1,0],\cc^4)$ we define a linear
operator:
$$L_\mu(\Phi)=A\Phi(0)+B_1\Phi(-\tau_1)+B_2\Phi(-\tau_2)$$
where A,$B_1$, $B_2$ are given by (7), (8) and a nonlinear
operator:
$$F(\mu, \Phi)=(0, -b_{12}\Phi_2(0)\Phi_1(0), \ds\f{1}{2}\rho_2\Phi_2^2(-\tau_1)+\ds\f{1}{6}\rho_3\Phi_2^3(-\tau_1),
0)^T+O(|\Phi|^4)$$ where $\Phi=(\Phi_1, \Phi_2, \Phi_3,
\Phi_4)^T$, $\rho_2=f^{''}(y_{10})$, $\rho_3=f^{'''}(y_{10}).$

By the Riesz representation theorem, there exists a matrix whose
components are bounded variation functions $\eta(\theta,\mu)$ with
$\theta\in[-\tau_{10},0]$ such that
$$L_\mu\Phi=\int_{-\tau_{10}}^0d\eta(\theta,\mu)\Phi(\theta),
\quad \theta\in[-\tau_{10},0].$$

We can choose
$$\eta(\theta,\mu)=\left\{\begin{array}{ll} \vspace{0.2cm}
A, & \th=0\\
B_2\delta(\theta+\tau_2), &
\th\in[-\tau_2,0)\\
B_1\delta(\theta+\tau_1), &
\th\in[-\tau_{10},-\tau_2).\end{array}\right.$$

For $\Phi\in C^1$ we define:
$$\ca(\mu)\Phi(\th)=\left\{\begin{array}{ll} \vspace{0.2cm}
\ds\f{d\Phi(\th)}{d\th}, & \th\in[-\tau_{10},0)\\
\int_{-\tau_{10}}^0d\eta(t,\mu)\Phi(t), &
\th=0,\end{array}\right.$$
$$R(\mu)\Phi=\left\{\begin{array}{ll} \vspace{0.2cm}
0, & \th\in[-\tau_{10},0)\\
F(\mu,\theta), & \th=0.\end{array}\right.$$

Then, we can rewrite (5) in the following vector form
$$\dot u_t=\ca(\mu)u_t+R(\mu)u_t\eqno(21)$$ where
$u_t=u(t+\theta)$, for $\theta\in[-\tau_{10},0]$.

As in [4] the bifurcating periodic solutions $u(t,\mu)$ of (21)
are indexed by a small parameters $\varepsilon$,
$\varepsilon\geq0$. The solution $u(t,\mu(\varepsilon))$ has
amplitude $O(\varepsilon)$, period $T(\varepsilon)$ and nonzero
Floquet exponent $\beta(\varepsilon)$ with $\beta(0)=0$, where
under our conditions $\mu$, T and $\beta$ have convergent
expansions:

$$\begin{array}{l}
\mu=\mu_2\varepsilon^2+\mu_4\varepsilon^4+\dots\\
T=\ds\f{2\pi}{\om_0}(1+\tau_2\varepsilon^2+\tau_4\varepsilon^4+\dots)\\
\beta=\beta_2\varepsilon^2+\beta_4\varepsilon^4+\dots
\end{array}$$

For $\Psi\in C^1([0,\tau_{10}], \cc^{*4})$, the adjoint operator
$\ca^*$ of $\ca$ is defined as:
$$\ca^*\Psi(s)=\left\{\begin{array}{ll} \vspace{0.2cm}
-\ds\f{d\Psi(s)}{ds}, & s\in(0, \tau_{10}]\\
\int^0_{-\tau_{10}}d\eta^T(t,0)\Psi(-t), &
s=0.\end{array}\right.$$

For $\Phi\in C([-\tau_{10}, 0], \cc^{4})$ and $\Psi\in
C^1([0,\tau_{10}], \cc^{*4})$ we define the following bilinear
form:
$$<\Psi(s),
\Phi(\th)>=\bar
\Psi(0)^T\Phi(0)-\int_{-\tau_{10}}^0\int_{\xi=0}^\theta\bar\Psi^T(\xi-\theta)d\eta(\theta)\Phi(\xi)d\xi,\eqno(22)$$
where $\eta(\theta)=\eta(\theta,0)$.

Then, it can verified that $\ca^*$ and $\ca$ are adjoint operators
with respect to this bilinear form.

For system (21) we have:

{\bf Proposition 6}. {\it If $\lambda_1=i\om_0$,
$\lam_2=\bar\lam_1$ then:

(i)The eigenvector of $\ca(0)$ corresponding to $\lam_1$ is
$$h(\th)=ve^{\lam_1\th},\quad \th\in[-\tau_{10}, 0]$$
where $v=(v_1, v_2, v_3, v_4)^T$,
$$v_1=0,v_2=-(\lam_1+d_2)(\lam_1+c_2),
v_3=-\rho e^{\lam_2\tau_{10}}(\lam_1+d_2), v_4=-\rho
e^{\lam_2\tau_{10}},$$ $\tau_0=\tau_{10}+\tau_{20}.$

(ii)The eigenvector of $\ca^*$ corresponding to $\lam_2$ is
$$h^{\ast}(s)=we^{\lam_1s},\quad s\in[0,\infty)$$ where
$w=(w_1, w_2, w_3, w_4)^T$,

$$\begin{array}{l}w_1=\eta, w_2=\ds\f{a_2+\lam_2}{b_1}\eta,
w_3=-\ds\f{e^{\lam_1\tau_2}b_{12}y_{10}(a_2+\lam_2)}{(c_2+\lam_2)(d_2+\lam_2)b_1}\eta,\\
w_4=-\ds\f{b_{12}y_{10}(a_2+\lam_2)e^{\lam_1\tau_2}}{b_1(d_2+\lam_2)(c_2+d_2)}\eta\end{array}$$
$$\begin{array}{l}\eta=\ds\f{a_2+\lam_2}{b_1}\bar v_2-(\bar
v_3-\rho\tau_{10}e^{\lam_1\tau_{10}}\bar
v_2)\ds\f{e^{\lam_1\tau_2}b_{12}y_{10}(a_2+\lam_2)}{(c_2+\lam_2)(d_2+\lam_2)b_1}-\\
-(\bar v_4-\tau_2e^{\lam_1\tau_2}\bar
v_3)\ds\f{b_{12}y_{10}(a_2+\lam_2)}{b_1(d_2+\lam_2)}\end{array}$$

(iii)With respect to (22) we have:
$$<h^{\ast}, h>=1,\quad <h^{\ast},\bar h>=<\bar h^{\ast}, h>=0, \quad <\bar h^{\ast}, \bar h>=1.$$}

Using the approach in [2], we next compute the coordinates to
describe the center manifold $\Omega_0$ at $\mu=0$. Let
$u_t=u(t+\th), \th\in[-\tau_{10},0)$, be the solution of system
(21) when $\mu=0$.

We define $$z(t)=<h^{\ast}, u_t>, \quad
w(t,\th)=u_t(\th)-2Re(z(t)h(\th)).$$

On the center manifold $\Omega_0$, we have:
$$w(t,\th)=w(z(t), \bar z(t), \th)$$ where
$$w(z,\bar z, \th)=w_{20}(\th)\ds\f{z^2}{2}+w_{11}(\th)z\bar
z+w_{02}(\th)\ds\f{\bar z^2}{2}+w_{30}(\th)\ds\f{z^3}{6}+\dots$$
in which $z$ and $\bar z$ are local coordinates for the center
manifold $\Omega_0$ in the direction of $h^{\ast}$ and $\bar
h^{\ast}$ and $w_{02}(\th)=\bar w_{20}(\th)$.

For solution $u_t\in\Omega_0$ of equation (21), as long as
$\mu=0$, we have:
$$\begin{array}{ll}
\dot z(t)& =\lam_1z(t)+\bar h^{\ast}(0)F(w(z(t),\bar
z(t),0)+2Re(z(t)h(0)))=\\
& \lam_1z(t)+g(z, \bar z)\end{array}$$ where
$$g(z, \bar z)=g_{20}\ds\f{z^2}{2}+g_{11}z\bar z+g_{02}\ds\f{\bar
z^2}{2}+g_{21}\ds\f{z^2\bar z}{2}+\dots$$

{\bf Proposition 7.} {\it For the system (21) we have:

(i)$$\begin{array}{l}g_{20}=-2b_{12}v_2v_4\bar w_2+\rho_2v_2^2\bar
w_3e^{2\lam_2\tau_{10}},\\
g_{11}=-b_{12}(v_2\bar v_4+\bar v_2v_4)\bar w_2+\rho_2v_2\bar v_2\bar w_3,\\
g_{02}=-2b_{12}\bar v_2\bar v_4\bar w_2+\rho_2\bar v_2^2\bar
w_3e^{2\lam_1\tau_{10}},\end{array}\eqno(23)$$

(ii)$$\begin{array}{l}
w_{20}(\th)=-\ds\f{g_{20}}{\lam_1}ve^{-\lam_1\th}-\ds\f{\bar
g_{02}}{3\lam_1}\bar ve^{\lam_2\th}+E_1e^{2\lam_1\th}\\
w_{11}(\th)=\ds\f{g_{11}}{\lam_1}ve^{\lam_1\th}-\ds\f{\bar
g_{11}}{\lam_1}ve^{\lam_2\th}+E_2,\\
\end{array}$$ where $E_1=(E_{11}, E_{21}, E_{31},
E_{41})^T$ and $E_2=(E_{12}, E_{22}, E_{32}, E_{42})^T$
$$\begin{array}{l}
E_{11}=0,
E_{21}=-\ds\f{\rho_2v_2^2}{\rho}+\ds\f{2\lam_1+c_2}{\rho}e^{-2\lam_2\tau_0}E_{41}\\
E_{31}=(2\lam_1+d_2)e^{-2\lam_2\tau_2}E_{41}\\
E_{41}=\ds\f{\rho_2v_2^2(2\lam_1+b_2+b_{12}y_{20})-2b_{12}v_2v_4\rho}{\rho
b_2y_{10}+(2\lam_1+b_2+b_{12}y_{20})(2\lam_1+c_2)e^{-2\lam_2\tau_0}}\end{array}$$
$$\begin{array}{l}
E_{12}=0,
E_{22}=-\ds\f{\rho_2v_2\bar v_2}{\rho}+\ds\f{c_2d_2}{\rho}E_2^4, E_{32}=d_2E_2^4\\
E_{42}=\ds\f{\rho b_{12}(v_2\bar v_4+\bar v_2v_4)-\rho_2v_2\bar
v_2(b_2+b_{12}y_{20})}{(b_2+b_{12}y_{20})c_2d_2+\rho
b_2y_{10}.}\end{array}$$

(iii)$$\begin{array}{ll}&g_{21} = -3b_{12}(\bar
v_2w_{420}(0)+2v_2w_{411}(0)+\bar
v_4w_{220}(0)+2w_{211}(0)v_4)\bar w_2+\\
& + \bar w_3[6\rho_2(2v_2e^{\lam_2\tau_1}-w_{211}(-\tau_1)+6\bar
v_2e^{\lam_1\tau_1}w_{220}(-\tau_1))+3\rho_3v^2_2e^{2\lam_2\tau_1}\bar
v_2e^{\lam_1\tau_1}],\\
\end{array}\eqno(24)$$ with $w_{20}(\theta)$=$(w_{120}(\theta)$, $w_{220}(\theta)$, $w_{320}(\theta)$, $w_{420}(\theta))^T$
and $w_{11}(\theta)$=$(w_{111}(\theta)$, $w_{211}(\theta)$,
$w_{311}(\theta)$, $w_{411}(\theta))^T$}

Based on the above analysis and calculation, we can see that each
$g_{ij}$ in (23), (24) are determined by the parameters and delay
from system (21). Thus, we can explicitly compute the following
quantities:
$$\begin{array}{l}
C_1(0)=\ds\f{i}{2\om_0}(g_{20}g_{11}-2|g_{11}|^2-\ds\f{1}{3}|g_{02}|^2)+\ds\f{g_{21}}{2}\\
\vspace{0.2cm} \mu_2=-\ds\f{Re(C_1(0))}{Re\lambda'(0)},
T_2=-\ds\f{Im(C_1(0))+\mu_2Im\lambda'(0)}{\om_0},
\beta_2=2Re(C_1(0)),\end{array}\eqno(25)$$ where $\lambda'(0)$ is
given by

$$\lambda'(0)=\left(\ds\f{r}{e^{\lam\tau}(3\lam^2+2p_2\lam+p_1-r)}\right)_{\lam=i\om_0,\tau=\tau_0}.$$

In summary, this leads to the following result:

\vspace{2mm} {\bf Theorem 1.} {\it In formulas (25), $\mu_2$
determines the direction of the Hopf bifurcation: if $\mu_2>0
(<0)$, then the Hopf bifurcation is supercritical (subcritical)
and the bifurcating periodic solutions exist for $\tau>\tau_{0}
(<\tau_{0})$; $\beta_2$ determines the stability of the
bifurcating periodic solutions: the solutions are orbitally stable
(unstable) if $\beta_2<0 (>0)$; and $T_2$ determines the period of
the bifurcating periodic solutions: the period increases
(decreases) if $T_2>0 (<0)$.}

\medskip

{\bf 3.2. The case $k_1(s)=\delta(s-\tau_1)$,
$k_2(s)=q_2e^{-q_2s}$, $\tau_1\geq0$, $q_2>0$.}

\hspace{0.6cm} For  $k_1(s)=\delta(s-\tau_1)$,
$k_2(s)=q_2e^{-q_2s}$, $\tau_1\geq0$, $q_2>0$, system (5) is given
by:
$$\begin{array}{l}
\vspace{0.1cm}
\dot u{}_1(t)=-a_2u_1(t),\\
\vspace{0.1cm}
\dot u{}_2(t)=b_1u_1(t)-(b_2+b_{12}y_{20})u_2(t)-b_{12}y_{10}u_4(t)-b_{12}u_2(t)u_4(t),\\
\dot u{}_3(t)=f(u_2(t-\tau_1)+y_{10}) -c_2(u_3(t)+x_{20}),\\
\dot u{}_4(t)=u_5(t)-d_2u_4(t),\\
\dot u{}_5(t)=q_2(u_3(t)-u_5(t)) .\end{array}\eqno(26)$$

We expand it in a Taylor series around the origin and neglect the
terms of higher order than the first order for the functions from
the right side of (26). We obtain:
$$\dot U(t)=A_{12}U(t)+B_{12}U(t-\tau_1),$$
where
$$A_{12}\!\!=\!\!\left(\!\!\!\!\begin{array}{ccccc}
\vspace{0.2cm}
-a_2 & 0 & 0 & 0 & 0\\
\vspace{0.2cm}
b_1 & -(b_2\!\!+\!\!b_{12}y_{20}) & 0 & -b_{12}y_{10} & 0\\
0 & 0 & -c_2 & 0 & 0\\
\vspace{0.2cm} 0 & 0 & 0 & -d_2 & 1\\
0 & 0 & q_2 & 0 & -q_2\\
\end{array}\right)\eqno(27)$$
$$B_{12}=\left(\begin{array}{ccccc} \vspace{0.2cm}
0 & 0 & 0 & 0 & 0\\
\vspace{0.2cm}
0 & 0 & 0 & 0 & 0\\
\vspace{0.2cm}
0 & \rho & 0 & 0 & 0\\
\vspace{0.2cm} 0 & 0 & 0 & 0 & 0\\
\vspace{0.2cm} 0 & 0 & 0 & 0 & 0\end{array}\right),\eqno(28) $$
with $U(t)=(u_1(t), u_2(t), u_3(t), u_4(t), u_5(t))$,
$U(t-\tau_1)=(u_1(t-\tau_1), u_2(t-\tau_1), u_3(t-\tau_1),
u_4(t-\tau_1), u_5(t-\tau_1))^T.$

Let $\tau_{10}^{\ast}$ given by (19) and
$\tau_{1}=\tau_{10}^{\ast}+\mu$, $\mu\in\rr$. We regard $\mu$ as
the bifurcation parameter.

For $\Phi\in C^1=C^1([-\tau_1,0],\cc^5)$ we define a linear
operator:
$$L_{12\mu}(\Phi)=A_{12}\Phi(0)+B_{12}\Phi(-\tau_1)$$
where $A_{12}$ ,$B_{12}$ are given by (27), (28) and a nonlinear
operator:
$$F_{12}(\mu, \Phi)=(0, -b_{12}\Phi_2(0)\Phi_1(0), \ds\f{1}{2}\rho_2\Phi_2^2(-\tau_1)+\ds\f{1}{6}\rho_3\Phi_2^3(-\tau_1),
0, 0)^T+O(|u|^4)$$ where $\Phi=(\Phi_1, \Phi_2, \Phi_3, \Phi_4,
\Phi_5)^T$, $\rho_2=f^{''}(y_{10})$, $\rho_3=f^{'''}(y_{10}).$

By the Riesz representation theorem, there exists a matrix whose
components are bounded variation functions $\eta(\theta,\mu)$ with
$\theta\in[-\tau_{10},0]$ such that
$$L_{12\mu}\Phi=\int_{-\tau_{10}^{\ast}}^0d\eta(\theta,\mu)\Phi(\theta),
\quad \theta\in[-\tau_{10}^{\ast},0].$$

We can choose
$$\eta_{12}(\theta,\mu)=\left\{\begin{array}{ll} \vspace{0.2cm}
A_{12}, & \th=0\\
B_{12}\delta(\theta+\tau_1), &
\th\in[-\tau_{10}^{\ast},0).\end{array}\right.$$

For $\Phi\in C^1$ we define:
$$\ca_{12}(\mu)\Phi=\left\{\begin{array}{ll} \vspace{0.2cm}
\ds\f{d\Phi(\th)}{d\th}, & \th\in[-\tau_{10}^{\ast},0)\\
\int_{-\tau_{10}^{\ast}}^0d\eta_{12}(t,\mu)\Phi(t), &
\th=0,\end{array}\right.$$
$$R_{12}(\mu)\Phi=\left\{\begin{array}{ll} \vspace{0.2cm}
0, & \th\in[-\tau_{10}^{\ast},0)\\
F_{12}(\mu,\theta), & \th=0.\end{array}\right.$$

Then, we can rewrite (26) in the following vector form
$$\dot u_t=\ca_{12}(\mu)u_t+R_{12}(\mu)u_t\eqno(29)$$ where
$u_t=u(t+\theta)$, for $\theta\in[-\tau_{10}^{\ast},0]$.

For $\Psi\in C^1([0,\tau_{10}^{\ast}], \cc^{*5})$, the adjoint
operator $\ca_{12}^*$ of $\ca$ is defined as:
$$\ca_{12}^*\Psi(s)=\left\{\begin{array}{ll} \vspace{0.2cm}
-\ds\f{d\Psi(s)}{ds}, & s\in(0, \tau_{10}^{\ast}]\\
\int^0_{-\tau_{10}^{\ast}}d\eta^T(t,0)\Psi(-t), &
s=0.\end{array}\right.$$

For $\Phi\in C([-\tau_{10}^{\ast}, 0], \cc^{5})$ and $\Psi\in
C^1([0,\tau_{10}^{\ast}], \cc^{*5})$ we define the following
bilinear form:
$$<\Psi(s),
\Phi(\th)>=\bar
\Psi(0)^T\Phi(0)-\int_{-\tau_{10}^{\ast}}^0\int_{\xi=0}^\theta\bar\Psi^T(\xi-\theta)d\eta_{12}(\theta)\Phi(\xi)d\xi,\eqno(30)$$
where $\eta(\theta)=\eta(\theta,0)$.

Then, it can verified that $\ca_{12}^*$ and $\ca_{12}$ are adjoint
operators with respect to this bilinear form.

For system (29) we have:

{\bf Proposition 8}. {\it If $\lambda_1=i\om_{10}$,
$\lam_2=\bar\lam_1$ then:

(i)The eigenvector of $\ca_{12}(0)$ corresponding to $\lam_1$ is
$$h(\th)=ve^{\lam_1\th},\quad \th\in[-\tau_{10}^{\ast}, 0]$$
where $v=(v_1, v_2, v_3, v_4, v_5)^T$,
$$v_1\!=\!0,v_2\!=\!\ds\f{(\lam_1\!+\!q_2)(\lam_1\!+\!c_2)}{\rho_1}e^{\lam_1\tau_1},
v_3=\lam_1\!+\!q_2, v_4=\!\ds\f{q_2}{\lam_1\!+\!d_2}, v_5\!=q_2.$$

(ii)The eigenvector of $\ca_{12}^*$ corresponding to $\lam_2$ is
$$h^{\ast}(s)=we^{\lam_1s},\quad s\in[0,\infty)$$ where
$w=(w_1, w_2, w_3, w_4, w_5)^T$,

$$\begin{array}{l}w_1=\ds\f{b_1}{(\lam_2+a_2)\eta}, w_2=\ds\f{1}{\eta},
w_3=-\ds\f{q_2b_{12}y_{10}}{(c_2+\lam_2)(d_2+\lam_2)(q_2+\lam_2)\eta},\\
w_4=-\ds\f{b_{12}y_{10}}{(d_2+\lam_2)\eta},
w_5=-\ds\f{b_{12}y_{10}}{(d_2+\lam_2)(q_2+\lam_2)\eta}\end{array}$$
$$\begin{array}{l}\eta\!=\!\bar v_2\!-\!\ds\f{q_2b_{12}y_{10}}{(c_2+\lam_2)(d_2+\lam_2)(q_2+\lam_2)}\bar
v_3\!-\!\ds\f{b_{12}y_{10}}{(d_2+\lam_2)}\bar
v_4\!-\!\ds\f{b_{12}y_{10}}{(d_2+\lam_2)(q_2+\lam_2)}\bar
v_5\-\!\\
-\ds\f{\rho_1q_2b_{12}y_{10}}{(c_2+\lam_2)(d_2+\lam_2)(q_2+\lam_2)\lam_2^2}
(e^{\lam_1\tau_{10}^{\ast}}-\tau_{10}^{\ast}\lam_2e^{\lam_1\tau_{10}^{\ast}}-1)\bar
v_2.\end{array}$$

(iii)With respect to (20) we have:
$$<h^{\ast}, h>=1,\quad <h^{\ast},\bar h>=<\bar h^{\ast}, h>=0, \quad <\bar h^{\ast}, \bar h>=1.$$}

Using the approach in [2], we next compute the coordinates to
describe the center manifold $\Omega_0$ at $\mu=0$. Let
$u_t=u(t+\th), \th\in[-\tau_{10}^{\ast},0)$, be the solution of
system (29) when $\mu=0$.

We define $$z(t)=<h^{\ast}, u_t>, \quad
w(t,\th)=u_t(\th)-2Re(z(t)h(\th)).$$

On the center manifold $\Omega_0$, we have:
$$w(t,\th)=w(z(t), \bar z(t), \th)$$ where
$$w(z,\bar z, \th)=w_{20}(\th)\ds\f{z^2}{2}+w_{11}(\th)z\bar
z+w_{02}(\th)\ds\f{\bar z^2}{2}+w_{30}(\th)\ds\f{z^3}{6}+\dots$$
in which $z$ and $\bar z$ are local coordinates for the center
manifold $\Omega_0$ in the direction of $h^{\ast}$ and $\bar
h^{\ast}$ and $w_{02}(\th)=\bar w_{20}(\th)$.

For solution $u_t\in\Omega_0$ of equation (29), as long as
$\mu=0$, we have:
$$\begin{array}{ll}
\dot z(t)& =\lam_1z(t)+\bar h^{\ast}(0)F(w(z(t),\bar
z(t),0)+2Re(z(t)h(0)))=\\
& \lam_1z(t)+g(z, \bar z)\end{array}$$ where
$$g(z, \bar z)=g_{20}\ds\f{z^2}{2}+g_{11}z\bar z+g_{02}\ds\f{\bar
z^2}{2}+g_{21}\ds\f{z^2\bar z}{2}+\dots$$

{\bf Proposition 9.} {\it For the system (29) we have:

(i)$$\begin{array}{l}g_{20}=-2b_{12}v_2v_4\bar w_2+\rho_2v_2^2\bar
w_3e^{2\lam_2\tau_{10}^{\ast}},\\
g_{11}=-b_{12}(v_2\bar v_4+\bar v_2v_4)\bar w_2+\rho_2v_2\bar v_2\bar w_3,\\
g_{02}=-2b_{12}\bar v_2\bar v_4\bar w_2+\rho_2\bar v_2^2\bar
w_3e^{2\lam_1\tau_{10}^{\ast}},\end{array}\eqno(31)$$

(ii)$$\begin{array}{l}
w_{20}(\th)=-\ds\f{g_{20}}{\lam_1}ve^{-\lam_1\th}-\ds\f{\bar
g_{02}}{3\lam_1}\bar ve^{\lam_2\th}+E_1e^{2\lam_1\th}\\
w_{11}(\th)=\ds\f{g_{11}}{\lam_1}ve^{\lam_1\th}-\ds\f{\bar
g_{11}}{\lam_1}ve^{\lam_2\th}+E_2,\\
\end{array}$$ where $E_1=(E_{11}, E_{21}, E_{31},
E_{41}, E_{51})^T$ and $E_2=(E_{12}, E_{22}, E_{32}, E_{42},
E_{52})^T$
$$\begin{array}{l}
E_{11}=0,
E_{21}=\ds\f{a_{22}F_{220}-a_{12}F_{320}}{a_{11}a_{22}+a_{12}a_{21}}, E_{31}=\ds\f{q_2+2\lam_1}{q_2}E_{51},\\
 E_{41}=\ds\f{1}{d_2+2\lam_1}E_{51}, E_{51}=\ds\f{a_{21}F_{220}+a_{11}F_{320}}{a_{11}a_{22}+a_{12}a_{21}}\end{array}$$
$$\begin{array}{l}
a_{11}=2\lam_1+b_2+b_{12}y_{10},
\quad a_{12}=\ds\f{b_{12}y_{10}}{d_2+2\lam_1},\\
a_{21}=\rho_1e^{2\lam_2\tau_{10}^{\ast}}, \quad
a_{22}=\ds\f{(c_2+2\lam_1)(q_2+2\lam_1)}{q_2},\end{array}$$

$$\begin{array}{l}
F_{220}=-2b_{12}v_2v_4,
F_{320}=\rho_2v_2^2e^{2\lam_2\tau_{10}^{\ast}},\end{array}$$

$$\begin{array}{l}
E_{12}=0,
E_{22}=\ds\f{c_{22}F_{211}-c_{12}F_{311}}{c_{11}c_{22}+c_{12}c_{21}}, E_{32}=E_{52}, E_{42}=\ds\f{1}{d_2}E_{52},\\
 E_{52}=\ds\f{c_{21}F_{211}+c_{11}F_{311}}{c_{11}c_{22}+c_{12}c_{21}}\end{array}$$

$$\begin{array}{l}
c_{11}=b_2+b_{12}y_{20}, c_{12}=b_{12}y_{10}, c_{21}=\rho_1,
c_{22}=c_2\\
F_{211}=-b_{12}(v_2\bar v_4+\bar v_2v_4), F_{311}=v_2\bar
v_2\rho_2.\end{array}$$

(iii)$$\begin{array}{ll}&g_{21} = -3b_{12}(\bar
v_2w_{420}(0)+2v_2w_{411}(0)+\bar
v_4w_{220}(0)+2w_{211}(0)v_4)\bar w_2+\\
& \!+\! \bar
w_3[6\rho_2(2v_2e^{\lam_2\tau_{10}^{\ast}}\!-\!w_{211}(\!-\!\tau_{10}^{\ast})\!+\!6\bar
v_2e^{\lam_1\tau_{10}^{\ast}}w_{220}(\!-\!\tau_{10}^{\ast}))\!+\!3\rho_3v^2_2e^{2\lam_2\tau_{10}^{\ast}}\bar
v_2e^{\lam_1\tau_{10}^{\ast}}],\\
\end{array}\eqno(32)$$ with $w_{20}(\theta)=(w_{120}(\theta)$, $w_{220}(\theta)$, $w_{320}(\theta)$, $w_{420}(\theta)$, $w_{520}(\theta))$
and $w_{11}(\theta)$=$(w_{111}(\theta)$, $w_{211}(\theta)$,
$w_{311}(\theta)$, $w_{411}(\theta)$, $w_{511}(\theta))$,
$\theta\in[-\tau_{10}^{\ast},0].$}

Based on the above analysis and calculation, we can see that each
$g_{ij}$ in (31), (32) are determined by the parameters and delay
from system (26). Thus, we can explicitly compute the following
quantities:
$$\begin{array}{l}
C_1(0)=\ds\f{i}{2\om_{10}}(g_{20}g_{11}-2|g_{11}|^2-\ds\f{1}{3}|g_{02}|^2)+\ds\f{g_{21}}{2}\\
\vspace{0.2cm} \mu_2=-\ds\f{Re(C_1(0))}{Re\lambda'(0)},
T_2=-\ds\f{Im(C_1(0))+\mu_2Im\lambda'(0)}{\om_{10}},
\beta_2=2Re(C_1(0)),\end{array}\eqno(33)$$ where $\lambda'(0)$ is
given by

$$\lambda'(0)=\left( \ds\f{\lambda
rq_2}{e^{\lambda\tau_1}(4\lambda^3\!+\!3(p_2\!+\!q_2)\lam^2\!+\!2(p_1\!+
\!q_2p_2)\lambda\!+\!p_0\!+\!p_1q_2)\!-\!rq_2\tau_1}\right)_{\lam\!=\!i\om_{10},\tau_1\!=\!\tau_{10}^{\ast}}.$$

We have:

\vspace{2mm} {\bf Theorem 2.} {\it In formulas (33), $\mu_2$
determines the direction of the Hopf bifurcation: if $\mu_2>0
(<0)$, then the Hopf bifurcation is supercritical (subcritical)
and the bifurcating periodic solutions exist for
$\tau_1>\tau_{10}^{\ast} (<\tau_{10}^{\ast})$; $\beta_2$
determines the stability of the bifurcating periodic solutions:
the solutions are orbitally stable (unstable) if $\beta_2<0 (>0)$;
and $T_2$ determines the period of the bifurcating periodic
solutions: the period increases (decreases) if $T_2>0 (<0)$.}

\section*{\normalsize\bf 4. Numerical examples.}

\hspace{0.6cm}For the numerical simulations we use Maple 9.5. In
this section, we consider system (6) with $a_1=2, a_2=0.55$,
$b_1=1$, $b_2=0.8$, $c_2=0.1$, $b_{12}=1.5$, $d_2=0.1$, $a=4$,
$n=2$. We obtain: $x_{10}= 3.636363636$, $y_{10}= 0.8347719895$,
$y_{20}\!=\! 2.370744013$, $x_{20}\!=\! 0.2370744013$.

In the first case, $k_1(s)=\delta(s-\tau_1)$,
$k_2(s)=\delta(s-\tau_2)$, for $\tau_2=3$, we have: $\omega_{10} =
0.1324013896$, $\mu_2\!=\! -0.4204703301$, $\beta_2\!=\!
0.2799153884$, $T_2\!=\! 0.0005051758260$,
$\tau_{0}\!=9.541873607\!$. Then the Hopf bifurcation is
subcritical and the bifurcating periodic solutions exist for
$\tau>\tau_{0}$; the solutions are orbitally unstable and the
period of the solution is increasing. The waveforms are displayed
in Fig1 and Fig2 and the phase plane diagrams of the state
variables $y_1(t)$, $y_2(t)$ and  $y_1(t-\tau)$, $y_1(t)$ are
displayed in Fig3 and Fig4:

\begin{center}
{\small \begin{tabular}{c|c} \hline Fig.1. $(t,y_1(t))$&Fig.2.
$(t,y_2(t))$\\
 \cline{1-2} \epsfxsize=5cm

\epsfysize=6cm

\epsffile{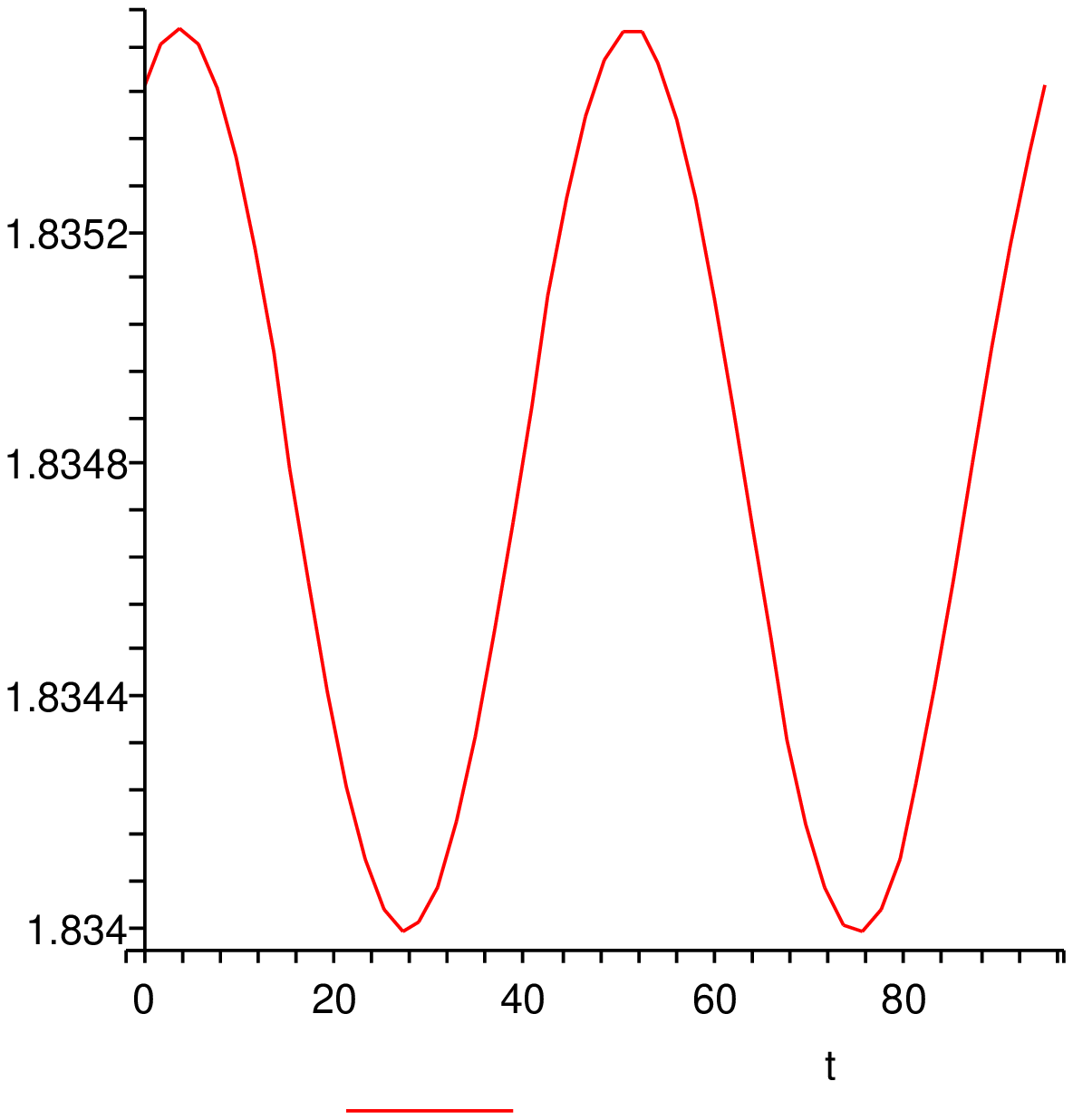} &

\epsfxsize=5cm

\epsfysize=6cm

\epsffile{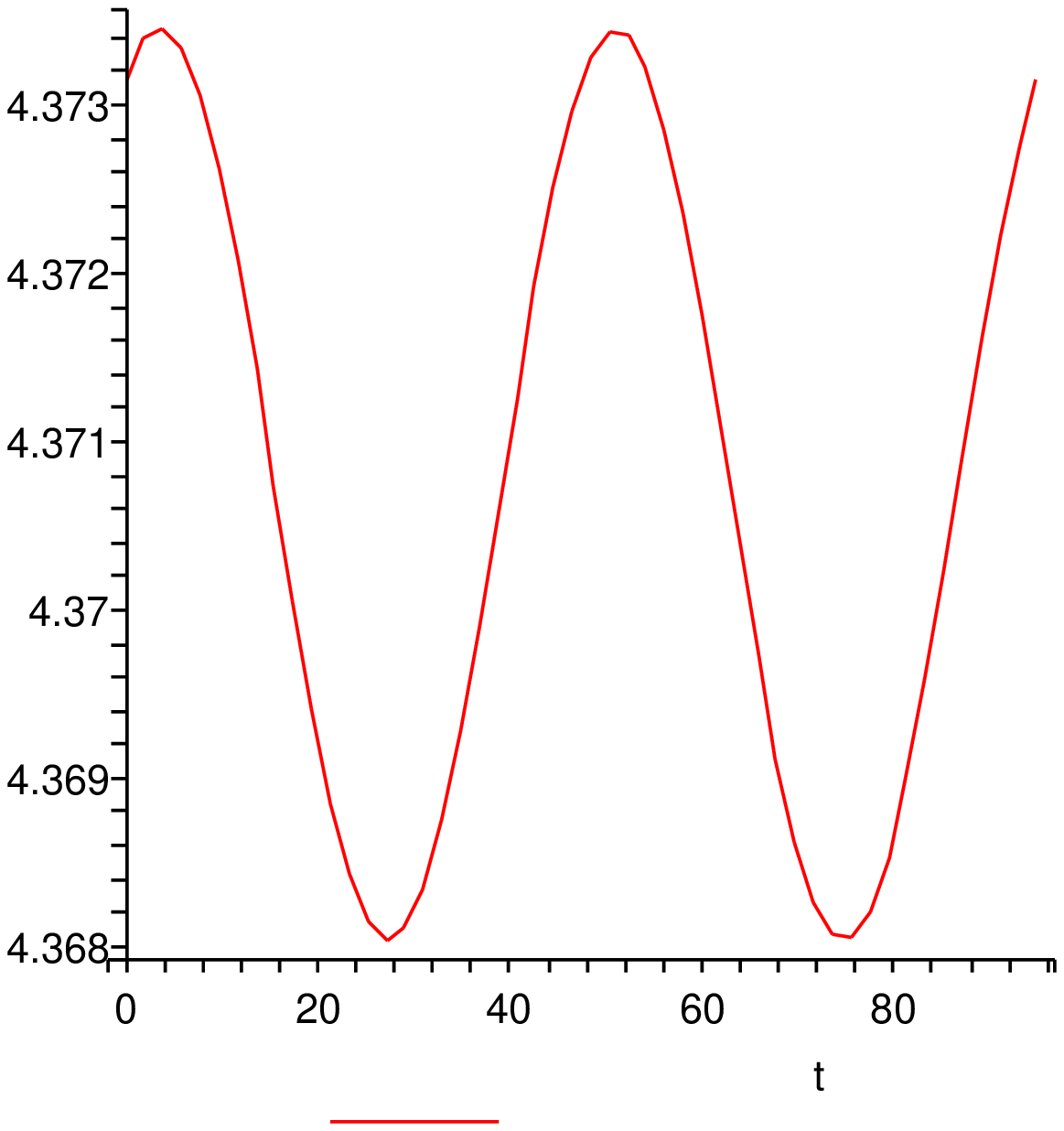}

\\
 \hline
\end{tabular}}
\end{center}

\medskip

\begin{center}
{\small \begin{tabular}{c|c} \hline Fig.3.
$(y_1(t),y_2(t))$&Fig.4.
 $(y_1(t-\tau),y_1(t))$\\
 \cline{1-2}

 \epsfxsize=5cm

\epsfysize=6cm

\epsffile{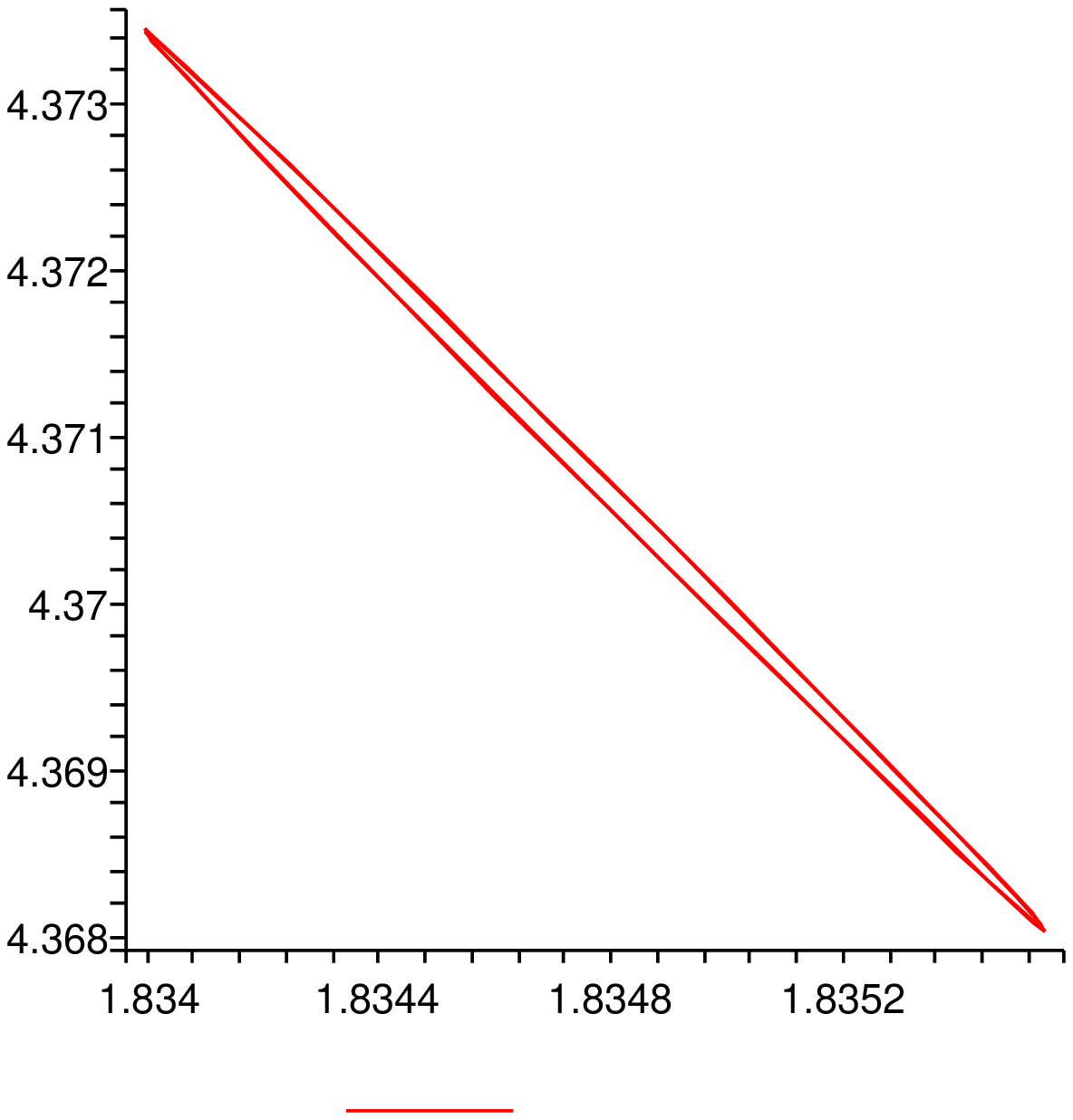} &

\epsfxsize=5cm

\epsfysize=6cm

\epsffile{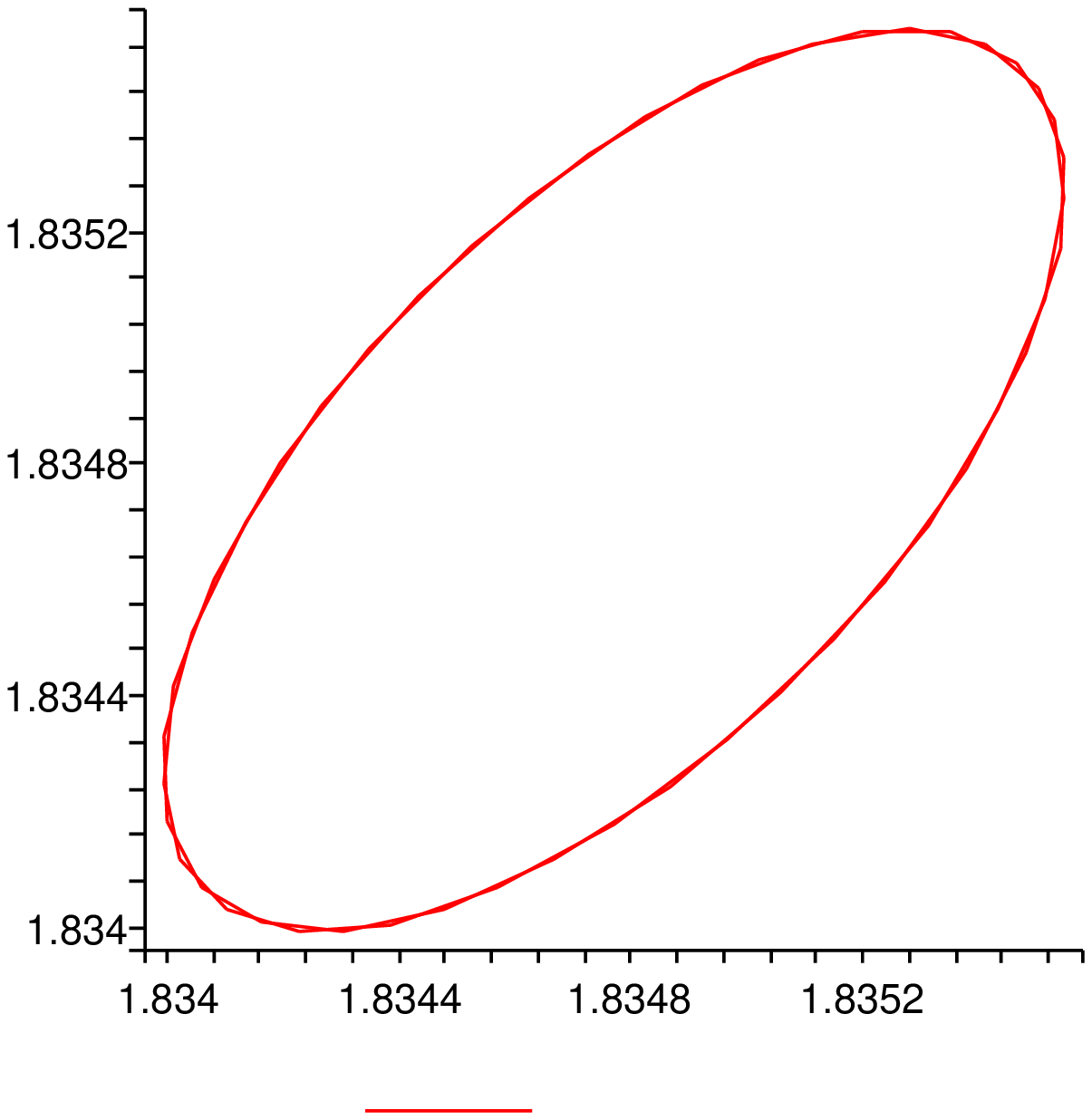}

\\
 \hline
\end{tabular}}
\end{center}

\medskip

In the second case, $k_1(s)=\delta(s-\tau_1)$,
$k_2(s)=q_2e^{-q_2s}$ for $q_2=0.5$, we have:  $\omega_{10} =
0.1290621026$, $\mu_2\!=\!-0.5993860816$, $\beta_2\!=\!
-0.7476750590$, $T_2\!=\! 0.1798944390$,
$\tau_{10}^{\ast}\!=32.37014890$. Then the Hopf bifurcation is
subcritical and the bifurcating periodic solutions exist for
$\tau_1>\tau_{10}^{\ast}$; the solutions are orbitally stable and
the period of the solution is increasing. The waveforms are
displayed in Fig5 and Fig6 and the phase plane diagrams of the
state variables $y_1(t)$, $y_2(t)$ and $y_1(t-\tau)$, $y_1(t)$ are
displayed in Fig7 and Fig8:

\begin{center}
{\small \begin{tabular}{c|c} \hline Fig.5. $(t,y_1(t))$&Fig.6.
$(t,y_2(t))$\\
 \cline{1-2} \epsfxsize=5cm

\epsfysize=6cm

\epsffile{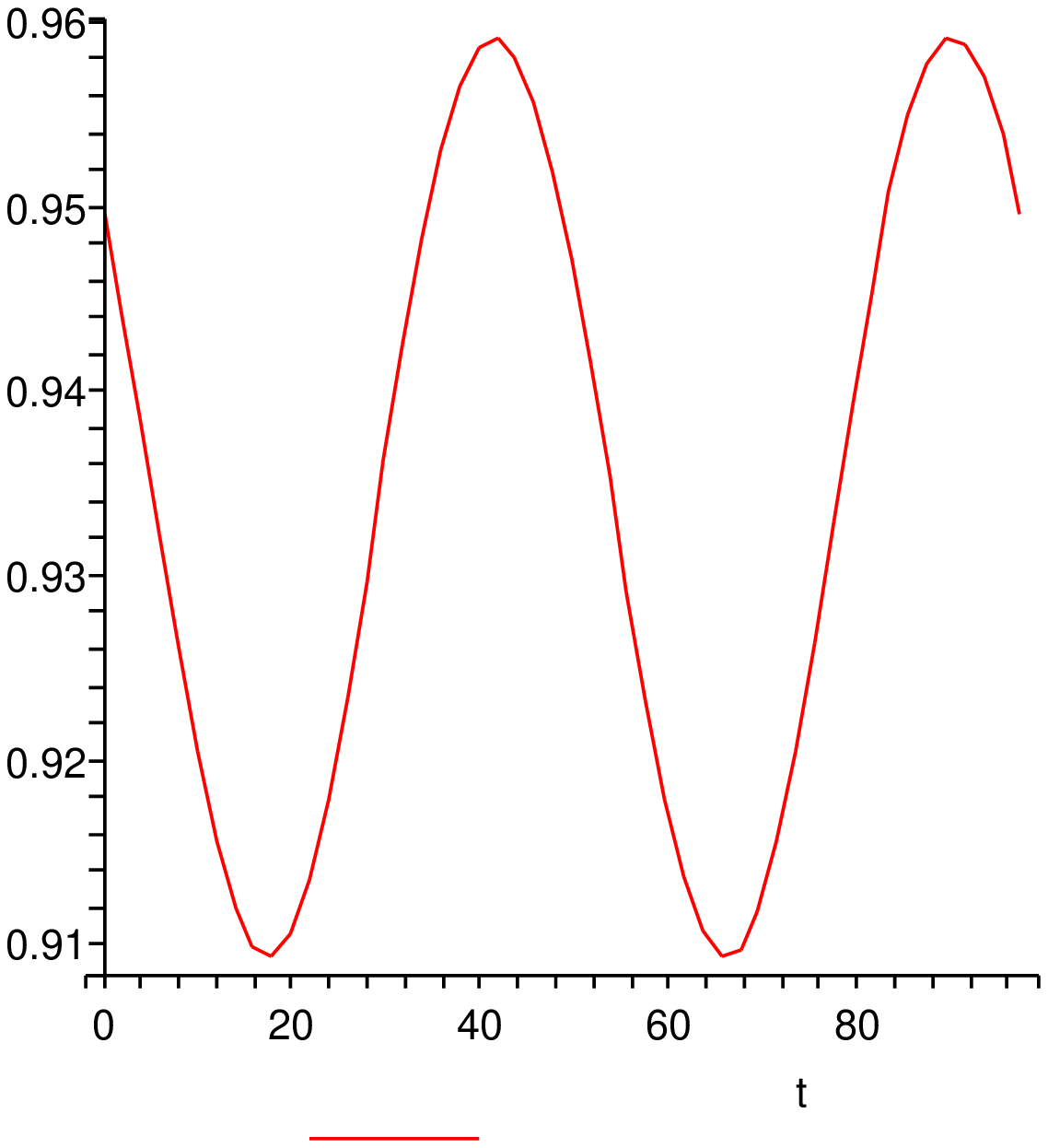} &

\epsfxsize=5cm

\epsfysize=6cm

\epsffile{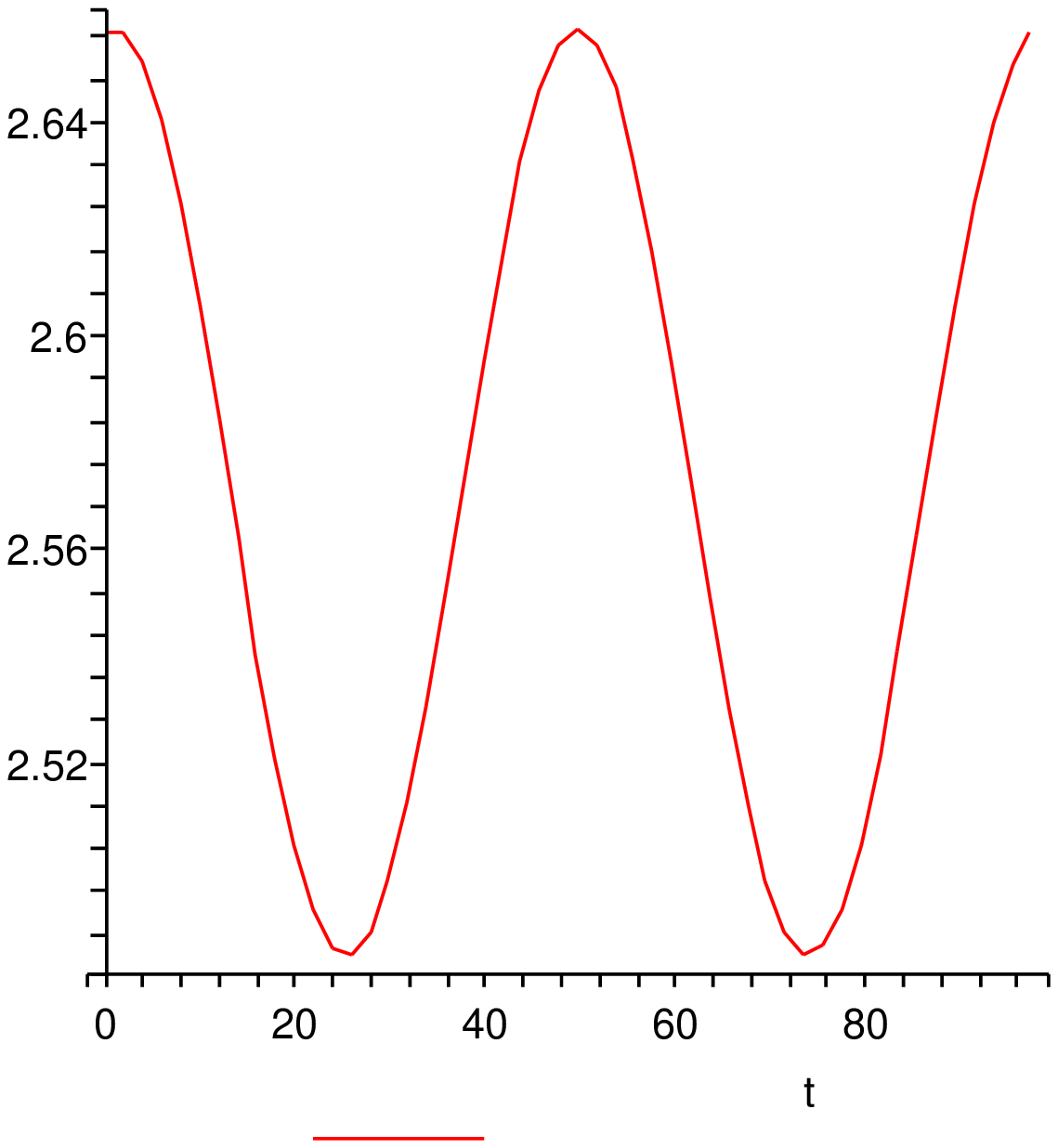}

\\
 \hline
\end{tabular}}
\end{center}

\medskip

\begin{center}
{\small \begin{tabular}{c|c} \hline Fig.7.
$(y_1(t),y_2(t))$&Fig.8.
 $(y_1(t-\tau),y_1(t))$\\
 \cline{1-2}

 \epsfxsize=5cm

\epsfysize=6cm

\epsffile{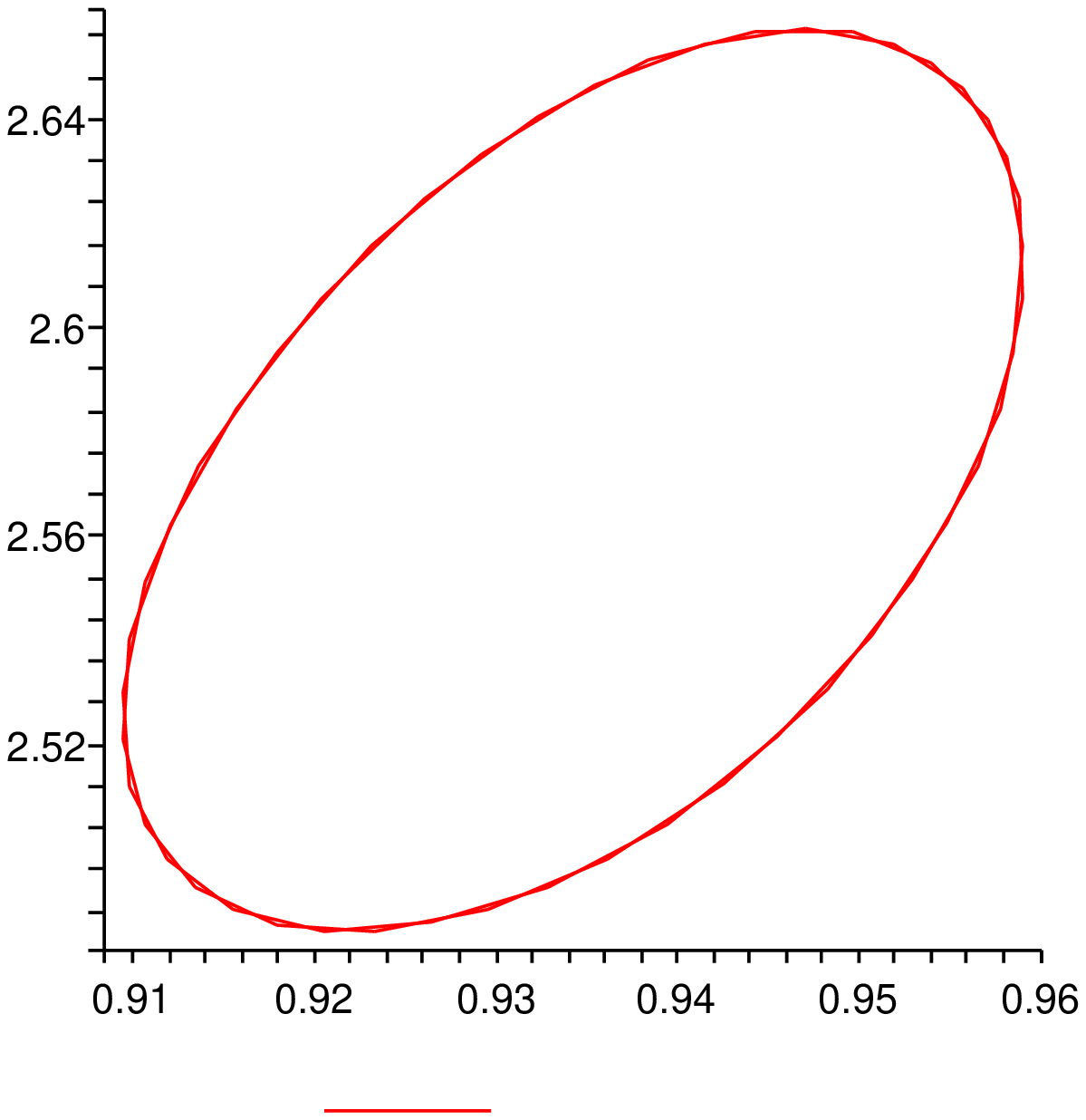} &

\epsfxsize=5cm

\epsfysize=6cm

\epsffile{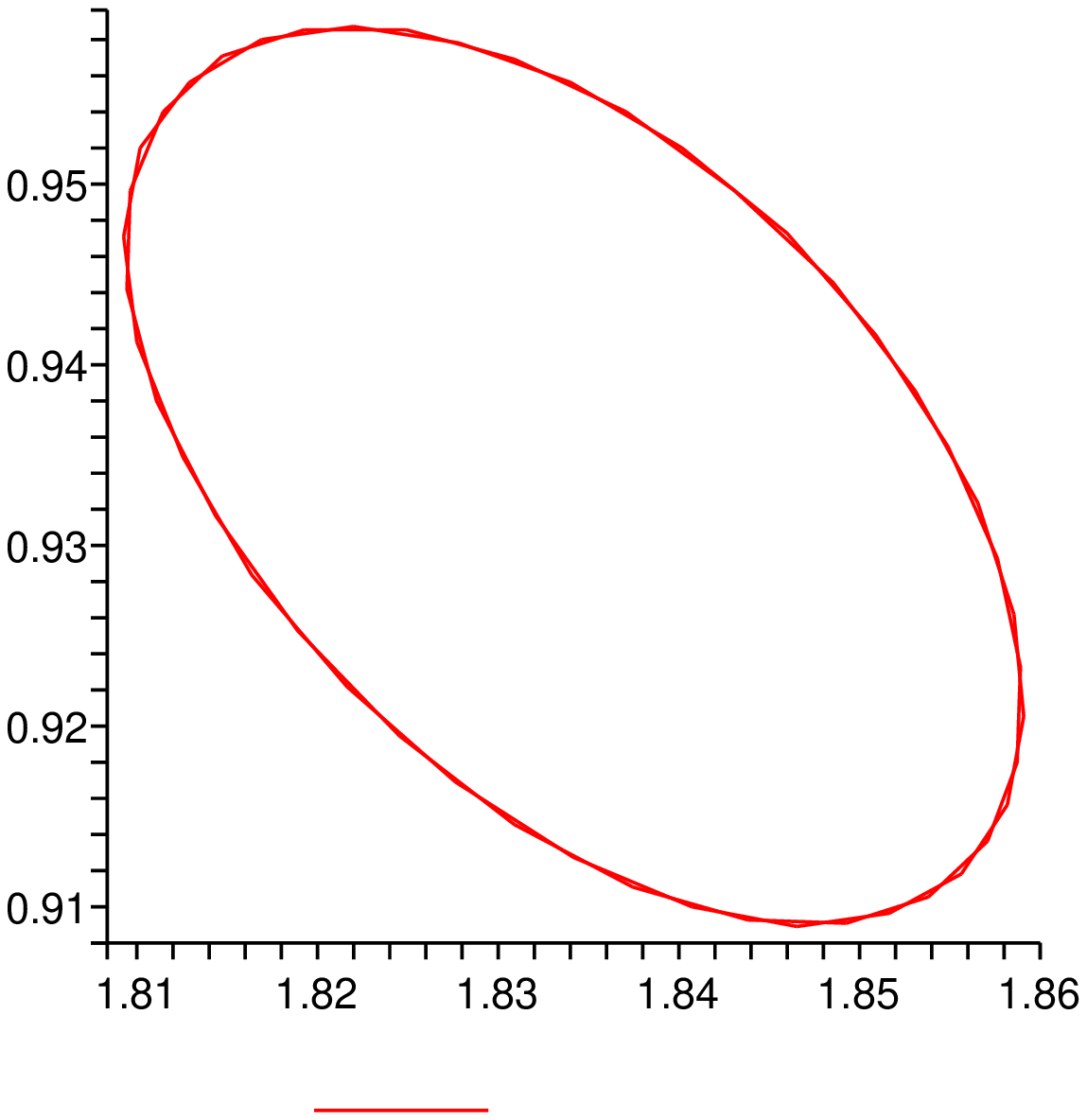}

\\
 \hline
\end{tabular}}
\end{center}

\medskip
\medskip



\section*{\normalsize\bf 5. Conclusions.}

\hspace{0.6cm} As in our previous models [11,13], we obtain an
oscillatory behavior similar to that observed experimentally [5].
The conclusion is not surprising, but is useful as this model
provides a more accurate approach of the interaction P53-Mdm2.

The improvements of the model from [11] done in the present paper,
proved usefulness as we obtained a smoother modelling of the
phenomenon and the oscillating behavior remained which as similar
with that from [5].

Using the method from this paper, we will do a qualitative
analysis of the model from [9] in our future papers.

\subsection*{Acknowledgements} This work is supported by Grant no 53(2006), "Ma\-the\-ma\-ti\-cal
models in hematology and P53-MDM2 dynamics with distributed time",
The National University Research Council from Ministry of
Education and Research of Romania.

\end{document}